\documentclass[11pt,reqno]{amsart}

\usepackage{amsmath, amsfonts, amsthm, amssymb, graphicx}

\theoremstyle{plain}
\newtheorem{theorem}{Theorem}[section]

\theoremstyle{definition}

\theoremstyle{remark}
\newtheorem{remark}{Remark}[section]

\numberwithin{equation}{section}

\newcommand{\e}{\varepsilon}

\renewcommand{\div}{\text{div}}
\renewcommand{\O}{\Omega}
\newcommand{\G}{\Gamma}
\renewcommand{\k}{\kappa}

\renewcommand{\d}{\partial}

\newcommand{\R}{{\mathbb R}}
\newcommand{\M}{{\mathbb M}}

\def\be{\begin{equation}}
\def\ee{\end{equation}}
\def\bes{\begin{equation*}}
\def\ees{\end{equation*}}
\def\bc{\begin{cases}}
\def\ec{\end{cases}}

\begin{document}

\title{Weak Continuity of the Gauss-Codazzi-Ricci System for Isometric Embedding}

\author{Gui-Qiang Chen \and  Marshall Slemrod \and Dehua Wang}
\address{G.-Q. Chen, School of Mathematical Sciences, Fudan University,
 Shanghai 200433, China; Department of Mathematics, Northwestern University,
         Evanston, IL 60208, USA.}
\email{\tt gqchen@math.northwestern.edu}

\address{M. Slemrod, Department of Mathematics, University of Wisconsin,
Madison, WI 53706, USA.} \email{\tt slemrod@math.wisc.edu}

\address{D. Wang, Department of Mathematics, University of Pittsburgh,
                Pittsburgh, PA 15260, USA.}
\email{\tt dwang@math.pitt.edu}

\date{\today}

%\thanks{}

\begin{abstract}
We establish the weak continuity of the Gauss-Coddazi-Ricci system
for isometric embedding with respect to the uniform $L^p$-bounded
solution sequence for $p>2$, which implies that the weak limit of
the isometric embeddings of the manifold is still an isometric
embedding. More generally, we establish a compensated compactness
framework for the Gauss-Codazzi-Ricci system in differential
geometry. That is, given
any sequence of approximate solutions to this system which is
uniformly bounded in $L^2$ and has reasonable bounds on the errors
made in the approximation (the errors are confined in a compact
subset of $H^{-1}_{\text{loc}}$), then the approximating sequence
has a weakly convergent subsequence whose limit is a solution of the
Gauss-Codazzi-Ricci system. Furthermore, a minimizing problem is
proposed as a selection criterion. For these, no restriction on the
Riemann curvature tensor is made.
\end{abstract}

\maketitle

\section{Introduction}

\medskip
The Gauss-Codazzi-Ricci system is a fundamental system of nonlinear
partial differential equations in differential geometry (cf.
\cite{BGY,BS,Eisenhart,Goenner,Greene,NM,Spivak}).
For example, the fundamental theorem of the surface theory indicates
that the existence of a local or global solution of the
Gauss-Codazzi-Ricci system can yield a local or global higher
dimensional isometric embedding. Therefore, it is important to
understand the behavior of this nonlinear system for solving
isometric embedding problems and other important geometric problems.
In general, the Gauss-Codazzi-Ricci system has no type, neither
purely hyperbolic nor purely elliptic.

\medskip
We are concerned with the weak continuity of the Gauss-Codazzi-Ricci
system and related compensated compactness framework for approximate
solutions to this system. In Chen-Slemrod-Wang \cite{CSW2}, we noted
that the Gauss-Codazzi equations for isometric embedding of $\M^2$
into $\R^3$ fall naturally within the formation of compensated
compactness. In this paper, we first show that this is also true in
the general case for the Gauss-Codazzi-Ricci system. One of our main
observations here is that the Codazzi and Ricci equations naturally
have the Div-Curl structure. Based on this observation, we establish
the week continuity of this system with respect to the uniform
$L^p$-bounded solution sequence for $p>2$, which implies that the
weak limit of the isometric embeddings of the manifold is still an
isometric embedding. This is reminiscent of the weak continuity of
determinants which plays an essential role in the theory of
polyconvexity  by Ball \cite{Ball} in nonlinear elasticity (also see
Dacorogna \cite{Dacorogna}, Evans \cite{Evans}, Morrey \cite{Mo66},
and M\"{u}ller \cite{Muller}). More generally, we establish a
stronger compensated compactness framework for the
Gauss-Codazzi-Ricci system. That is, given
any sequence of approximate solutions to this system which is
uniformly bounded in $L^2$, and has reasonable bounds on the errors
made in the approximation (the errors are confined in a compact
subset of $H^{-1}_{\text{loc}}$), then the approximating sequence
has a weakly convergent subsequence whose limit is still a solution
of the Gauss-Codazzi-Ricci system. For these, no restriction on the
Riemann curvature tensor is made.

\medskip
A long-standing fundamental problem in differential geometry is the
existence of
local (and if possible global) embeddings of a $d$-dimensional
Riemannian manifold $\M^d$, $d\ge 3$, into the Euclidean space
$\R^N$ with optimal dimension $N$. As noted in
Han-Hong \cite{HanHong}, the first global existence of smooth
embeddings was given by Nash \cite{Nash1}, but the best result as of
this time is the following theorem of G\"unther \cite{Gunther}: {\em
Any smooth $d$-dimensional compact Riemannian manifold admits a
smooth (i.e. $C^\infty$) isometric embedding in $\R^N$ for
$N=\frac{1}{2}\max\{d(d+5), d(d+3)+10\}$}.
Needless to say, it is of considerable interest to know if
G\"unther's dimension $N$ is optimal. In a similar vein, we could
try to formulate a selection or ``admissibility" criterion to choose
one of the possibly infinite embeddings provided by G\"unther's
theorem. Within the realm of surface theory and elastic manifolds,
this has been recently considered in \cite{Witten1, Witten2} where
the selection is done by minimizing an integral of norm of the
second fundamental form. Indeed, this seems a natural approach for
selection in the general case and is even in the same spirit of
Dafermos's entropy rate criterion \cite{Dafermos-book}. In Section
4, we propose a minimizing problem  as a selection criterion and
show by the compensated compactness framework that any minimizing
sequence has a subsequence in $L^p, p>2$, which converges weakly to
a minimizer that satisfies the Gauss-Codazzi-Ricci system. Since any
sequence of isometric embeddings of $\M^d$ into $\R^N$ (say given by
G\"unther's theorem) must satisfy the equations exactly, this
implies that the problem of minimizing the $L^p$-norms of the second
fundamental form and the connection form on the normal bundle
(sometimes called torsion coefficients \cite{BS})
does have a solution within the class of weak solutions of the
Gauss-Codazzi-Ricci system, hence yielding an isometric immersion of
$W^{2,p}$ class for $p>2$.

\bigskip
\section{The Gauss-Codazzi-Ricci System for Isometric Embedding of $\M^d$ into $\R^N$}

\medskip
In this section, we use the following conventional notation:
\begin{equation*}
\begin{split}
g_{ij}: &\quad\text{given metric of the Riemannian manifold},\\
\G_{ij}^k: &\quad\text{Christoffel symbols},\\
R_{ijkl}: &\quad\text{Riemann curvature tensor},\\
h_{ij}^a: & \quad\text{Coefficients of the second fundamental form},\\
\k_{lb}^a: &\quad\text{Coefficients of the connection form (torsion
coefficients) on the normal bundle},
\end{split}
\end{equation*}
where the indices $a, b, c$ run from $1$ to $N$, and $i, j, k, l, m,
n$ run from $1$ to $d\ge 3$.

\medskip
For given metric $g_{ij}$, the Christoffel symbols are
$$
\G_{ij}^{k}=\frac12g^{kl}\left(\d_j g_{il}+\d_i g_{jl}-\d_l
 g_{ij}\right),
$$
which depend on the first derivatives of $(g_{ij})$, and the Riemann
curvature tensor is
$$
R_{ijkl}=g_{lm}\left(\d_k\G^{m}_{ij}-\d_j\G^{m}_{ik}
+\G^{n}_{ij}\G^{m}_{nk}-\G^{n}_{ik}\G^{m}_{nj}\right),
$$
which depends on $(g_{ij})$ and its first and second derivatives,
where $(g^{kl})$ denotes the inverse of $(g_{ij})$ and
$\partial_j=\partial_{x_j}$. We denote $|g|=det(g_{ij})$.

\subsection{The Gauss-Codazzi-Ricci System}

As is well-known in Riemannian Geometry, the isometric embedding
problems for $d$-dimensional Riemannian manifolds into the Euclidean
space $\R^N$ can be reduced as the solvability problems of the
Gauss-Codazzi-Ricci system of nonlinear partial differential
equations with the following form:

\medskip
The Gauss equations:
\begin{equation} \label{G1}
h_{ji}^a h_{kl}^a-h_{ki}^a h_{jl}^a=R_{ijkl};
\end{equation}

The Codazzi equations: \be\label{C1} \frac{\d h_{lj}^a}{\d
x^k}-\frac{\d h_{kj}^a}{\d x^l} + \G_{lj}^mh_{km}^a-\G_{kj}^m
h_{lm}^a +\k_{kb}^a h_{lj}^b-\k_{lb}^ah_{kj}^b=0; \ee

The Ricci equations: \be\label{R1} \frac{\d\k_{lb}^a}{\d
x^k}-\frac{\d\k_{kb}^a}{\d x^l}
-g^{mn}\left(h^a_{ml}h^b_{kn}-h^a_{mk}h^b_{ln}\right)
+\k_{kc}^a\k_{lb}^c-\k_{lc}^a\k_{kb}^c=0.
\ee

\medskip
Notice that the coefficients of the second fundamental form are
symmetric: \be\label{s1} h_{ij}^a=h_{ji}^a,
\end{equation}
while the coefficients of the connection form on the normal bundle
are antisymmetric:
\begin{equation}\label{s2}
\k_{kb}^a=-\k_{ka}^b.
\end{equation}
In particular, the antisymmetry of $\kappa_{kb}^a$ implies
$$
\k_{ka}^a=-\k_{ka}^a,
$$
and so
$$
\k_{ka}^a=0.
$$
Thus, the $a$th column of the $d\times d$ matrix $\k^a$ is zero.

\medskip
When $d=3$, the Janet dimension $N=\frac{d(d+1)}{2}=6$ (cf. Janet
\cite{Janet}). Then
\begin{equation*}
\k^1= \begin{bmatrix}
           0 & \k_{12}^1 & \k_{13}^1\\
           0 & \k_{22}^1 & \k_{23}^1\\
           0 & \k_{32}^1 & \k_{33}^1
          \end{bmatrix}, \quad
\k^2= \begin{bmatrix}
           -\k_{12}^1 & 0 & \k_{13}^2 \\
           -\k_{22}^1 & 0 & \k_{23}^2 \\
           -\k_{32}^1 & 0 & \k_{33}^2
          \end{bmatrix}, \quad
\k^3= \begin{bmatrix}
           -\k_{13}^1 & -\k_{13}^2 & 0 \\
           -\k_{23}^1 & -\k_{23}^2 & 0 \\
           -\k_{33}^1 & -\k_{33}^2 & 0
           \end{bmatrix}.
\end{equation*}

\subsection{The Div-Curl Structure of the Codazzi and Ricci
Equations}

In this section we present one of our main observations on the
features of the Codazzi and Ricci equations: the Div-Curl structure,
which leads to the weak continuity of the system.

\smallskip
For $w=(w_1, w_2, \cdots, w_d)$,
$$
\text{curl}\, w:=(\partial_jw_i-\partial_i w_j)_{1\le i,j\le d}
$$
is a $d\times d$ matrix field.

\medskip
From the Codazzi equations \eqref{C1}, for $k<l$, they possess the
form:
$$
\frac{\d h^a_{lj}}{\d x^k}-\frac{\d h^a_{kj}}{\d x^l} + l.o.t=0,
$$
or
\begin{equation}\label{C11}
\div(\underbrace{\overbrace{0, \cdots, h_{lj}^a}^{k}, 0, \cdots,
-h_{kj}^a}_{l}, 0,\cdots, 0)+ l.o.t=0, \ee
and
\begin{equation}\label{C12}
\text{curl}(h^a_{1j}, h^a_{2j}, \cdots, h_{dj}^a)+ l.o.t=0,
\end{equation}
where $l.o.t$ represents the lower-order terms without involving
derivatives in the equation.

Similarly, we observe that the identical form of the Ricci equations
\eqref{R1} can also be written as
\begin{equation}\label{R11}
\div(\underbrace{\overbrace{0, \cdots, 0, \kappa_{lb}^a}^{k}, 0,
\cdots, -\kappa_{kb}^a}_l, 0, \cdots, 0)+ l.o.t=0,
\end{equation}
and
\begin{equation}\label{R12}
\text{curl}(\k^a_{1b}, \k^a_{2b}, \cdots, \k^a_{db})+ l.o.t=0.
\end{equation}

\medskip
Now replacing $a$ by $b$, and $j$ by $i$ in the Codazzi
equations \eqref{C11}--\eqref{C12}, we obtain
\begin{equation}\label{C13}
\div(\underbrace{\overbrace{0, \cdots, h_{li}^b}^{k}, 0, \cdots,
-h_{ki}^b}_{l}, 0,\cdots, 0)+ l.o.t=0, \ee
and
\begin{equation}\label{C14}
\text{curl}(h^b_{1i}, h^b_{2i}, \cdots, h^b_{di})+ l.o.t=0.
\end{equation}
Similarly, replacing $a$ by $b$ and $b$ by $c$ in the Ricci
equations \eqref{R11}--\eqref{R12}, we have
\begin{equation}\label{R13}
\div(\underbrace{\overbrace{0, \cdots, 0, \kappa_{lc}^b}^{k}, 0,
\cdots, -\kappa_{kc}^b}_l, 0, \cdots, 0)+ l.o.t=0,
\end{equation}
and
\begin{equation}\label{R14}
\text{curl}(\k^b_{1c}, \k^b_{2c}, \cdots, \k^b_{dc})+ l.o.t=0.
\end{equation}

\iffalse The div-curl lemma asserts that if we have an approximate
sequence $h^a_{ij}, \k^c_{lb}$ uniformly bounded in $L^p(\O),
\O\subset\R^d$, for which the right hand sides of our div, curl
equalities lie in compact subsets of $H^{-1}_{\text{loc}}(\O)$, then
scalar product of any pair of vectors in our Codazzi and Ricci list
is weakly continuous in $L^p(\O)$.
\fi

\medskip
One of our main observations is that the scalar product of the two
vector fields in the rewritten forms \eqref{C11}--\eqref{R14} yield
the nonlinear quantities in the lower-order terms in the
Gauss-Codazzi-Ricci system \eqref{G1}--\eqref{R1}: Forms \eqref{C11}
and \eqref{C14} yield
\begin{equation}\label{qq-1}
h_{lj}^ah_{ki}^b-h_{kj}^ah_{li}^b;
\end{equation}
forms \eqref{R11} and \eqref{R14} yield
\begin{equation}\label{qq-2}
\kappa_{kb}^a\kappa_{lc}^b-\kappa_{lb}^a\kappa_{kc}^b;
\end{equation}
and forms \eqref{R12} and \eqref{C13} yield
\begin{equation}\label{qq-3}
\kappa_{kb}^ah_{li}^b-\kappa_{lb}^a h_{ki}^b.
\end{equation}
This observation is essential for us to establish the weak
continuity of the Gauss-Codazzi-Ricci system in \S 3.

\bigskip
\section{Weak Continuity and Compensated Compactness Framework}

\medskip
In this section we establish the weak continuity of the
Gauss-Codazzi-Ricci system and related compensated compactness
framework for approximate solutions to the system via the Div-Curl
lemma (see Murat \cite{Murat1} and Tartar \cite{Tartar1}).

The Div-Curl lemma is a basic result in the compensated compactness
theory for the weak continuity of the scalar product of two vector
fields (cf. \cite{Dacorogna,Evans,Murat1,Murat2,Tartar1,Tartar2})
and is closely related with the Hodge decomposition.

\begin{theorem}[Div-Curl Lemma]
Let $\Omega\subset\R^d, d\ge 2,$ be open bounded. Let $p, q>1$ such
that $\frac{1}{p}+\frac{1}{q}=1$. Assume that, for any
$\varepsilon>0$, two fields $u^\varepsilon\in L^p(\Omega; \R^d)$ and
$v^\varepsilon\in L^q(\Omega; \R^d)$ satisfy the following:
\begin{enumerate}\renewcommand{\theenumi}{\roman{enumi}}
\item
\label{weakSol-def-i1} $u^\varepsilon\rightharpoonup u$ weakly in
$L^p(\Omega;\R^d)$ as $\varepsilon\to 0$;
\item
\label{weakSol-def-i2} $v^\varepsilon\rightharpoonup v$ weakly in
$L^q(\Omega;\R^d)$ as $\varepsilon\to 0$;
\item  ${\rm div}\, u^\varepsilon$ are confined in a
compact subset of $W^{-1, p}_{loc}(\Omega; \R)$;
\item  ${\rm curl}\, v^\varepsilon$ are confined in a
compact subset of $W^{-1, q}_{loc}(\Omega; \R^{d\times d})$.
\end{enumerate}
Then the scalar product of $u^\e$ and $v^\e$ are weakly continuous:
$$
u^\varepsilon\cdot v^\varepsilon\longrightarrow u\cdot v
$$
in the sense of distributions.
\end{theorem}

Based on our observation of the Div-Curl structure of the Codazzi
and Ricci equations, we employ the Div-Curl lemma to formulate the
following compensated compactness framework.

\medskip
Let a sequence of vector fields $(h_{ij}^{a,\varepsilon},
\kappa_{lb}^{a,\varepsilon})({\bf x})$, defined on an open bounded
subset $\Omega\subset \R^d$, satisfy the following Framework (A):

\medskip
{\rm (A.1)}\, $\|(h_{ij}^{a,\varepsilon},
\kappa_{lb}^{a,\varepsilon})\|_{L^2(\Omega)}\le C$ for some $C>0$
independent of $\varepsilon>0$;

\smallskip
{\rm (A.2)}\, $\frac{\d h_{lj}^{a,\e}}{\d x^k}-\frac{\d
h_{kj}^{a,\e}}{\d x^l}$ and $\frac{\d\k_{lb}^{a,\e}}{\d
x^k}-\frac{\d\k_{kb}^{a,\e}}{\d x^l}$ are confined in a compact set
in $H_{loc}^{-1}(\Omega)$;

\smallskip
{\rm (A.3)\, There exist $o^\varepsilon_j(1), j=1,2,3$, with
$o^\varepsilon_j(1)\to 0$ in the sense of distributions as
$\varepsilon\to 0$ such that
\begin{equation} \label{g1-b}
\begin{split}
&\frac{\d h_{lj}^{a,\e}}{\d x^k}-\frac{\d h_{kj}^{a,\e}}{\d x^l}
+\G_{lj}^mh_{km}^{a,\e}-\G_{kj}^mh_{lm}^{a,\e}
+\k_{kb}^{a,\e} h_{lj}^{b,\e}-\k_{lb}^{a,\e}h_{kj}^{b,\e}=o^{\varepsilon}_1(1),\\
&\frac{\d\k_{lb}^{a,\e}}{\d x^k}-\frac{\d\k_{kb}^{a,\e}}{\d x^l}
-g^{mn}\left(h^{a,\e}_{ml}h^{b,\e}_{kn}-h^{a,\e}_{mk}h^{b,\e}_{ln}\right)
+\k_{kc}^{a,\e}\k_{lb}^{c,\e}-\k_{lc}^{a,\e}\k_{kb}^{c,\e}=o^{\varepsilon}_2(1),
\end{split}
\end{equation}
and
\begin{equation}\label{g2-b}
h_{ji}^{a,\e} h_{kl}^{a,\e}-h_{ki}^{a,\e} h_{jl}^{a,\e}=R_{ijkl}
+o^{\varepsilon}_3(1).
\end{equation}

Then we have
\begin{theorem}[Compensated compactness framework] \label{T}
$\quad$ Let a sequence of vector fields $(h_{ij}^{a,\varepsilon},
\kappa_{lb}^{a,\varepsilon})$ satisfy Framework {\rm (A)}. Then
there exists a subsequence (still labeled) $(h_{ij}^{a,\varepsilon},
\kappa_{lb}^{a,\varepsilon})$ that converges weakly in $L^2(\Omega)$
to $(h_{ij}^{a}, \kappa_{lb}^{a})$ as $\varepsilon\to 0$ such that

\begin{enumerate}
\renewcommand{\theenumi}{\roman{enumi}}
\item  $\|(h_{ij}^{a},
\kappa_{lb}^{a})\|_{L^2(\Omega)}\le C$;

\item the quadratic terms in \eqref{G1}--\eqref{R1} are
weakly continuous with respect to the subsequence
$(h_{ij}^{a,\varepsilon}, \kappa_{lb}^{a,\varepsilon})$ that
converges  to $(h_{ij}^{a}, \kappa_{lb}^{a})$ weakly in
$L^2(\Omega)$ as $\varepsilon\to 0$;

\item the limit vector field $(h_{ij}^{a},\kappa_{lb}^{a})$ satisfies
the Gauss-Codazzi-Ricci system \eqref{G1}--\eqref{R1}.
\end{enumerate}

\noindent That is, the limit vector field $(h_{ij}^{a},
\kappa_{lb}^{a})$ is a weak solution to the Gauss-Codazzi-Ricci
system \eqref{G1}--\eqref{R1}.
\end{theorem}

\begin{proof}
By assumption (A.1), there exists a subsequence (still denoted)
$(h_{ij}^{a,\varepsilon}, \kappa_{lb}^{a,\varepsilon})$ and a vector
field $(h_{ij}^{a}, \kappa_{lb}^{a})\in L^2(\Omega)$ such that
\begin{equation}\label{3.3}
(h_{ij}^{a,\varepsilon}, \kappa_{lb}^{a,\varepsilon})
\rightharpoonup (h_{ij}^{a}, \kappa_{lb}^{a}) \qquad\mbox{in}\,\,
L^2(\Omega),
\end{equation}
and
\begin{equation}\label{3.4} \|(h_{ij}^{a},
\kappa_{lb}^{a})\|_{L^2(\Omega)}\le C.
\end{equation}

By the Div-Curl structure, observed in \S 2.2, assumption (A.2)
implies that
\begin{equation}\label{C111}
\div(\underbrace{\overbrace{0, \cdots, h_{lj}^{a,\e}}^{k}, 0,
\cdots, -h_{kj}^{a,\e}}_{l}, 0,\cdots, 0), \quad
\text{curl}(h^{a,\e}_{1j}, h^{a,\e}_{2j}, \cdots, h^{a,\e}_{dj})
\end{equation}
and
\begin{equation}\label{R111}
\div(\underbrace{\overbrace{0, \cdots, 0, \kappa_{lb}^{a,\e}}^{k},
0, \cdots, -\kappa_{kb}^{a,\e}}_l, 0, \cdots, 0), \quad
\text{curl}(\k^{a,\e}_{1b}, \k^{a,\e}_{2b}, \cdots, \k^{a,\e}_{db})
\end{equation}
are confined in a compact set in \, $H^{-1}_{loc}(\Omega)$.

\medskip
By exchanging the indices, we also have
\begin{equation}\label{C131}
\div(\underbrace{\overbrace{0, \cdots, h_{li}^{b,\e}}^{k}, 0,
\cdots, -h_{ki}^{b,\e}}_{l}, 0,\cdots, 0), \quad
\text{curl}(h^{b,\e}_{1i}, h^{b,\e}_{2i}, \cdots, h^{b,\e}_{di})
\end{equation}
and
\begin{equation}\label{R131}
\div(\underbrace{\overbrace{0, \cdots, 0, \kappa_{lc}^{b,\e}}^{k},
0, \cdots, -\kappa_{kc}^{b,\e}}_l, 0, \cdots, 0),\quad
\text{curl}(\k^{b,\e}_{1c}, \k^{b,\e}_{2c}, \cdots, \k^{b,\e}_{dc})
\end{equation}
are confined in a compact set in $\, H^{-1}_{loc}(\Omega)$.

\medskip
Using the Div-Curl lemma, Theorem 3.1, we conclude that the weak
continuity of the nonlinear quadratic quantities in the
Gauss-Codazzi-Ricci system with respect to the sequence
$(h_{ij}^{a,\varepsilon}, \kappa_{lb}^{a,\varepsilon})$:
\begin{eqnarray}
&&h_{lj}^{a,\e}h_{ki}^{b,\e}-h_{kj}^{a,\e}h_{li}^{b,\e}\,\,\rightharpoonup\,\,
h_{lj}^{a}h_{ki}^{b}-h_{kj}^{a}h_{li}^{b},\label{3.9}\\
&&\kappa_{kb}^{a,\e}\kappa_{lc}^{b,\e}-\kappa_{lb}^{a,\e}\kappa_{kc}^{b,\e}
\,\, \rightharpoonup\,\,
\kappa_{kb}^a\kappa_{lc}^b-\kappa_{lb}^a\kappa_{kc}^b,\label{3.10}\\
&&\kappa_{kb}^{a,\e}h_{li}^{b,\e}-\kappa_{lb}^{a,\e} h_{ki}^{b,\e}
\,\, \rightharpoonup\,\, \kappa_{kb}^ah_{li}^b-\kappa_{lb}^a
h_{ki}^b\label{3.11}
\end{eqnarray}
in the sense of distributions as $\varepsilon\to 0$.

\medskip

Combining \eqref{3.3}--\eqref{3.4} with \eqref{3.9}--\eqref{3.11},
we conclude that the weak limit vector field $(h_{ij}^{a},
\kappa_{lb}^{a})$ of the sequence $(h_{ij}^{a,\varepsilon},
\kappa_{lb}^{a,\varepsilon})$ satisfy the Gauss-Codazzi-Ricci system
\eqref{G1}--\eqref{R1} in the sense of distributions, that is, the
limit vector field $(h_{ij}^{a}, \kappa_{lb}^{a})$ is a weak
solution of \eqref{G1}--\eqref{R1}.
\end{proof}

As a corollary, we conclude the weak continuity of the
Gauss-Codazzi-Ricci system with respect to the uniform $L^p$-bounded
solution sequence for $p>2$.

\begin{theorem}[Weak Continuity]
Let $(h_{ij}^{a,\varepsilon}, \kappa_{lb}^{a,\varepsilon})$ be a
sequence of solutions to the Gauss-Codazzi-Ricci system
\eqref{G1}--\eqref{R1}, which is uniformly bounded in $L^p$, $p>2$.
Then the weak limit vector field $(h_{ij}^{a}, \kappa_{lb}^{a})$ of
the sequence $(h_{ij}^{a,\varepsilon}, \kappa_{lb}^{a,\varepsilon})$
in $L^p$ is still a solution to \eqref{G1}--\eqref{R1}.
\end{theorem}

\begin{proof}
Since the solution sequence $(h_{ij}^{a,\varepsilon},
\kappa_{lb}^{a,\varepsilon})$ is uniformly bounded in $L^p$, $p>2$:
\begin{equation}\label{p-bound} \|(h_{ij}^{a,\e},
\kappa_{lb}^{a,\e})\|_{L^p(\Omega)}\le C,
\end{equation}
for some $C>0$ independent of $\e$, then there exists  a subsequence
(still denoted) $(h_{ij}^{a,\varepsilon},
\kappa_{lb}^{a,\varepsilon})$ and a vector field $(h_{ij}^{a},
\kappa_{lb}^{a})\in L^p(\Omega)$ such that
$$
(h_{ij}^{a,\varepsilon}, \kappa_{lb}^{a,\varepsilon})\,\,
\rightharpoonup \,\, (h_{ij}^{a}, \kappa_{lb}^{a})
\qquad\mbox{in}\,\, L^p(\Omega),
$$
and
$$
\|(h_{ij}^{a}, \kappa_{lb}^{a})\|_{L^p(\Omega)}\le C.
$$

\medskip Then we conclude from \eqref{p-bound} that all the
lower-order terms for the solution sequence
$(h_{ij}^{a,\varepsilon}, \kappa_{lb}^{a,\varepsilon})$ in the
Gauss-Codazzi-Ricci system \eqref{G1}--\eqref{R1} are uniformly
bounded in $L^{p/2}, p>2$. This implies that
\begin{equation}\label{h-compact}
\frac{\d h_{lj}^a}{\d x^k}-\frac{\d h_{kj}^a}{\d x^l}, \,\,
\frac{\d\k_{lb}^a}{\d x^k}-\frac{\d\k_{kb}^a}{\d x^l} \quad
\mbox{are confined in a compact set in $H_{loc}^{-1}(\Omega)$}.
\end{equation}

\medskip Since the domain $\Omega\subset\R^d$ is bounded, the
uniform bound in \eqref{p-bound} implies the uniform bound of
$(h_{ij}^{a,\e}, \kappa_{lb}^{a,\e})$ in $L^2(\Omega)$.  By the
compensated compactness framework (Theorem 3.2), we conclude that
the limit vector field is a weak solution of \eqref{G1}--\eqref{R1},
which implies the weak continuity of the system.
\end{proof}

\begin{remark}
The weak continuity of the Gauss-Codazzi-Ricci system implies that,
for $p>2$, the weak limit of a sequence of isometric embeddings of
the $d$-dimensional manifold $\M^d$ into $\R^N$ as surfaces with
corresponding uniform $L^p$-bounded sequence
$(h_{ij}^{a,\varepsilon}, \kappa_{lb}^{a,\varepsilon})$ is still an
isometric embedding as a surface in $\R^N$. The requirement $p>2$ is
to ensure the $H^{-1}$-compactness in \eqref{h-compact} to deal with
the nonhomogeneous terms.
\end{remark}

\bigskip
\section{Minimization Problem}
\medskip

In this section, as an example, we show that the solution sequence
$(h_{ij}^{a,\varepsilon}, \kappa_{lb}^{a,\varepsilon})$ for the weak
continuity in Theorem 3.3 can be obtained from a selection
criterion.

\begin{theorem}
There exists a minimizer $(h_{ij}^a, \k_{lb}^a)$ for the
minimization problem: \be \label{M1} \min_S \|(h,\k)\|_{L^p(\O)}^p:=
\min_S
\int_{\O}\sqrt{|g|}\left((h_{ij}h_{ij})^\frac{p}2+(\k_{lb}\k_{lb})^\frac{p}2\right)
dx, \ee where $S$ is the set of weak solutions to the
Gauss-Codazzi-Ricci system \eqref{G1}--\eqref{R1}.
\end{theorem}

\begin{proof}
Clearly, $S$ is non-empty by G\"{u}nther's theorem in \cite{Gunther}
(also see the statement in \S 1 above). A minimizing sequence
provides the desired $L^p$-norm for the weak continuity theorem
(Theorem 3.3). Since the $L^p$-norm is convex, which is weakly lower
semicontinuous, any minimizing sequence has a subsequence in
$L^p(\O)$ that converges weakly to a minimizer which satisfies the
Gauss-Codazzi-Ricci system \eqref{G1}--\eqref{R1}.
\end{proof}

Notice that any sequence of isometric embeddings of $\M^d$ into
$\R^N$ as surfaces (say, given by G\"unther's theorem) must satisfy
the Gauss-Codazzi-Ricci equations \eqref{G1}--\eqref{R1}. This
implies that the problem of minimizing the $L^p$-norms of the second
fundamental form and the connection form on the normal bundle
does have a solution within the class of weak solutions of the
Gauss-Codazzi-Ricci system \eqref{G1}--\eqref{R1}, hence yielding an
isometric immersion of $W^{2,p}$ class for $p>2$ for $\M^d$ into
$\R^N$ as a surface.

\bigskip\bigskip\bigskip\bigskip\bigskip
{\bf Acknowledgments.}
 Gui-Qiang Chen's research was supported in
part by the National Science Foundation under Grants DMS-0807551,
DMS-0720925, and DMS-0505473, and the Natural Science Foundation of
China under Grant NSFC-10728101. Marshall Slemrod's research was
supported in part by
 the National Science Foundation under Grant DMS-0647554.
 Dehua Wang's research was supported in part by the National Science
Foundation under Grant DMS-0604362, and by the Office of Naval
Research under Grant N00014-07-1-0668. This paper was written as
part of the International Research Program on Nonlinear Partial
Differential Equations at the Centre for Advanced Study at the
Norwegian Academy of Science and Letters in Oslo during the academic
year 2008--09; It was finalized when the authors participated in the
SQuaRE on ``Isometric Embedding of Higher Dimensional Riemannian
Manifolds'', which was held at the American Institute of
Mathematics, Palo Alto, California, March 16--20, 2009.

\end{document}